\documentclass[12pt,reqno]{amsart}
    \usepackage{amssymb,amsmath,amsthm,newlfont}
    \usepackage[dvips,final]{graphicx}
\usepackage[T1]{fontenc}
\usepackage{a4wide}

\numberwithin{equation}{section}
\newcommand{\T}{\mathbb{T}}

\newcommand{\N}{\mathbb{N}}
\newcommand{\R}{\mathbb{R}}

\newcommand{\E}{\mathbb{E}}

\newcommand{\D}{\mathbb{D}}

\newcommand{\cD}{{\mathcal{D}}}

\newcommand{\cA}{{\mathcal{A}}}
\newcommand{\cH}{{\mathcal{H}}}

\newcommand{\cL}{{\Lambda_\omega}}

\newtheorem{thm}{\bf Theorem}[]

\newtheorem{lem}{\bf Lemma}

\begin{document}
\title{Outer functions in analytic weighted Lipschitz algebras}
\date{\today}
\subjclass[2000]{primary 46E20; secondary 30C85, 47A15.}
\author[B. Bouya]{Brahim Bouya}

\curraddr{Madretschstr.52, CH-2503, Biel.}

\thanks{A part of this work was partially supported by ANR FRAB: ANR-09-BLAN-0058-02}

\email{brahimbouya@gmail.com}
\maketitle

\begin{abstract}
{We give a complete description of outer functions in the analytic weighted Lipschitz algebras by their moduli in the boundary, with respect to any modulus of continuity.}
\end{abstract}

\section{\bf Introduction and statement of main result.}

Let $\D$ be the open unit disk of  the complex plane and $\T$ its boundary.
In all what follows, we let $h:\T\rightarrow \R^+$  be a nonnegative continuous function such that
\begin{eqnarray}\label{logint}
\int_{-\pi}^{\pi}\log h(\theta) d\theta>-\infty,
\end{eqnarray}
where $h(\theta):=h(e^{i\theta})$ for all real numbers $\theta\in\R$.
The outer function associated to $h$ is defined by
$$O_{h}(z):=\exp\{u_{_h}(z)+iv_{_h}(z)\},\qquad z\in\D,$$
where $$u_{_h}(z):=\displaystyle\frac{1}{2\pi}\int_{-\pi}^{\pi}\text{Re}\Big(\frac{e^{i\varphi}+z}{e^{i\varphi}-z}\Big)\log h(\varphi)d\varphi,  \quad z\in\D,$$
and $v_{_h}$ is the harmonic conjugate of the harmonic function $u_{_h}$ given by
$$v_{_h}(z):=\frac{1}{2\pi}\int_{-\pi}^{\pi}\text{Im}\Big(\frac{e^{i\varphi}+z}{e^{i\varphi}-z}\Big)\log h(\varphi)d\varphi, \qquad z\in\D,$$
see for instance the books \cite{Gar,Hof}.
We denote by $\cH^{\infty}(\D)$ the space of all bounded analytic functions in $\D.$
It is well known that $O_h$  is a function in $\cH^{\infty}(\D)$ which is not necessarily in $\cA(\D),$ the disk algebra of all functions in $\cH^{\infty}(\D)$ that are continuous up to the boundary. In Theorem \ref{continu} below we will
give a necessary and sufficient condition on $h$ so that $O_h\in\cA(\D).$
Since outer functions are completely constructed from their values of their moduli in the boundary, one may ask which properties
on $h$ ensure that $O_h$ belongs to a given space of analytic functions.
The famous Carleson's formula \cite{Car1} in the standard Dirichlet space $\cD,$ of analytic functions on $\D$ with a
square area integrable modulus of the derivatives, establishes in fact a necessary and sufficient condition only on $h$ ensuring $O_h\in\cD,$
see also \cite[Chapter 7.4]{EKMR}.
Later, in \cite[Theorem 3.1]{Shi2} an analogue of Carleson's formula was given by Shirokov  for some spaces of analytic functions smooth up to the boundary, concerning this result we can see also Theorem \ref{polycase} below.

We now let $\omega$ be  a {\it modulus of
continuity}, i.e., a nondecreasing continuous real-valued function on $[0,2]$ with $\omega(0)=0,$ $\omega(1)=1$ and such that $t\mapsto\omega(t)/t$ is nonincreasing. As examples for such $\omega,$ we can consider the following ones
\begin{equation*}
\chi_\alpha(t):=t^\alpha, \qquad 0\leq t\leq 2,
\end{equation*}
\begin{equation*}
\varphi_{_\beta}(t):=a_\beta\log^{-\beta}\big(\frac{2e^{\beta}}{t}\big), \qquad 0< t\leq 2,
\end{equation*}
and
\begin{equation*}
\psi_{_\beta}(t):=b_\beta\log^{-\beta} \log\big(\frac{2e^{1+\beta}}{t}\big), \qquad 0< t\leq 2,
\end{equation*}
where $0<\alpha\leq1$ and $\beta>0$ are fixed real numbers, and $a_\beta,b_\beta>0$ are the positive constants such that our regularity condition $\varphi_{_\beta}(1)=\psi_{_\beta}(1)=1$ holds.
Given a modulus of continuity $\omega,$ the analytic weighted Lipschitz algebra $\cL$ is the space of all functions $f\in\cA(\D)$ such that
$$\|f\|_{\omega,\D}:=\sup_{z\in\D}|f(z)|+ \sup_{z,w\in\D \atop z\neq w}\frac{|f(z)-f(w)|}{\omega(|z-w|)}<+\infty.$$
The above result of Shirokov has been shown only for modulus of continuity satisfying the following stronger condition
\begin{equation}\label{c}
\int_{0}^{s}\frac{\omega(t)}{t}dt+s\int_{s}^{2}\frac{\omega(t)}{t^2}dt\leq c\omega(s), \qquad s\in]0,2],
\end{equation}
where $c>0$ is a constant. However, this condition is not satisfied by modulus of continuity of slow growth such as for both $\varphi_{_\beta}$ and $\psi_{_\beta}.$

Our main result is Theorem \ref{asmahan1} below, in which we give necessary and sufficient conditions only on $h$ so that  $O_{h}\in\cL,$ with respect to any modulus of continuity.
We denote by $\cL(\T)$ the space of all continuous functions $g$ in $\T$ such that
$$\|g\|=\|g\|_{\omega,\T}:=\sup_{z\in\T}|g(z)|+ \sup_{z,w\in\T \atop z\neq w}\frac{|g(z)-g(w)|}{\omega(|z-w|)}<+\infty.$$
For two distinct points $e^{i\theta}$ and $e^{i\varphi}$ such that $-\pi\leq\theta<\varphi<\pi,$ we define the arc $\gamma:=(e^{i\theta},e^{i\varphi})$ of $\T$ joining the points  $e^{i\theta}$ and $e^{i\varphi}$  as follows $$\gamma:=\{e^{is}\ :\ \theta< s<\varphi\}.$$
Associated to a function $h\in\cL(\T)$ and a positive number $\eta\leq1,$ we define the following set of arcs
$$\Gamma_{h,\eta}:=\{\gamma\subseteq\T\ :\ 0<|\gamma|\leq\omega^{*}\big(\frac{\eta h(\gamma)}{2\|h\|}\big)\}$$
where  $|\gamma|$ designs the arc length of the arc $\gamma,$
$$h(\gamma):=\inf\{h(s)\ :\ e^{is}\in\gamma\},$$
and $\omega^*$ is the right inverse of $\omega:$
$$\omega^{*}(u):=\inf\{t\in[0,1]\ :\ \omega(t)=u\},\qquad u\in[0,1].$$
We set
$$A_h(\gamma):=\lim_{\varepsilon\rightarrow 0}\int_{e^{is}\in\gamma}\int_{\varepsilon}^{\omega^{*}(\frac{h(\gamma)}{2\|h\|})}
 \frac{h(s+t) h(s-t)-h^2(s)}{t^{2}}dsdt,\qquad \gamma\in\Gamma_{h},$$
where $\Gamma_{h}:=\Gamma_{h,1}.$
We also define the following Korenblum's function
$$a_{h}(\theta):=\int_{|\varphi|\leq \pi\atop h(\varphi)\leq \frac{1}{2}h(\theta)}\frac{\log\frac{h(\theta)}{h(\varphi)}}{|e^{i\varphi}-e^{i\theta}|^2}d\varphi,\qquad e^{i\theta}\in\T\setminus \E_h,$$
where $\E_h$ designs the zero set of $h.$

\begin{thm}\label{asmahan1}
Let $h:\T\mapsto \R^+$ be a continuous nonnegative function satisfying \eqref{logint} and such that $h\in\cL(\T)$ for some modulus of continuity $\omega.$ Let $\rho\geq 1$ be a real number.
The outer function $O^{\rho}_{h}$ belongs to $\Lambda_{\omega}$ if and only if $h$ satisfies the following two conditions
\begin{equation}\label{condition1}
\sup_{e^{i\theta}\in\T\setminus\E_h}\frac{h^{\rho}(\theta)}{\omega(\min\{1, 1/a_{h}(\theta)\})}<+\infty.
\end{equation}
There exists a positive number $\eta\leq 1$ such that
\begin{equation}\label{condition2}
 \sup_{\gamma\in\Gamma_{h^{\rho},\eta}}\frac{|A_{h^{\rho}}(\gamma)|}{h^{\rho}(\gamma)\omega(|\gamma|)}<+\infty.
 \end{equation}
\end{thm}

\begin{proof}
See section \ref{proftheorem}.
\end{proof}

The raison for which we are interested in Theorem \ref{asmahan1} to the situations when $\rho\geq 1$ is that because as much as $\rho$ is bigger the condition \eqref{condition1} gets more chance to be satisfied , see Theorem \ref{coro2} below. However, the condition \eqref{condition2} generally depends on the local comportment of $h$ on $\T\setminus\E_h,$ see Lemma \ref{smouth} below.

For the limiting situation of the algebras $\cL,$ we obtain the following theorem.

\begin{thm}\label{continu}
Let $h:\T\mapsto \R^+$ be a continuous nonnegative function satisfying \eqref{logint}. We have $O_h\in\cA(\D)$
if and only if for each point $\xi\in\T\setminus\E_h$ there exist a neighborhood $\mathcal{U}\subseteq\T$ of $\xi$ and a positive number $\lambda\leq\pi$
such that
\begin{equation}\label{continuity}
\lim_{|\gamma|\rightarrow0}\lim_{\varepsilon\rightarrow 0}\int_{e^{is}\in\gamma}\int_{\varepsilon}^{\lambda}
\frac{h(s+t) h(s-t)-h^2(s)}{t^{2}}dsdt=0,
\end{equation}
uniformly with respect to all arcs $\gamma$ of $\mathcal{U}.$
\end{thm}

\begin{proof}
See section \ref{sect4}.
\end{proof}

The remainder  of this paper is organized as follows: In the next section we give the statement of some interesting particular cases (Theorem \ref{polycase} and Theorem \ref{coro2}) of Theorem \ref{asmahan1}. Section \ref{apps} is devoted to presenting an application of Theorem \ref{coro2}. In sections \ref{sect4} and \ref{proofpoly} we will give respectively the proof of Theorem \ref{continu} and Theorem \ref{polycase}. In the last section we will prove Theorem \ref{asmahan1}.

\section{ \bf Special cases.}

A particular case of Theorem \ref{asmahan1} provide us with the next one, in which we encounter (assertion \eqref{first}) a result due to Shirokov \cite[Theorem 3.1]{Shi3}.

\begin{thm}\label{polycase}
Let $\omega$ be modulus of continuity satisfying \eqref{c}.
Let $h:\T\mapsto \R^+$ be a continuous nonnegative function satisfying \eqref{logint} and such that $h\in\cL(\T).$ Then
$ O_{h}\in\cL$ if and only if one of the following assertions holds
  \begin{eqnarray}\label{first}
\sup_{e^{i\theta}\in\T\setminus\E_{h}}\frac{h(\theta)}{\omega\big(\min\{1, a^{-1}_{h}(\theta)\}\big)}<+\infty.
\end{eqnarray}
There exists a positive number $\delta\leq 1$ such that
 \begin{eqnarray}\label{second}
 \sup_{z\in\D \atop 1-|z|= \omega^{*}(\frac{\delta h(z/|z|)}{2\|h\|})}\big|\log\frac{|O_{h}(z)|}{h(z/|z|)}\big|<+\infty.
 \end{eqnarray}

\end{thm}

\begin{proof}
See section \ref{proofpoly}.
\end{proof}

In the next section we will give an application to the following theorem.

\begin{thm} \label{coro2}
Let $\omega$ be a modulus of continuity such that
\begin{equation}\label{cond}
\inf_{0< t\leq2}\frac{\omega(t^{2})}{\omega^{\rho}(t)}>0,
\end{equation}
for some number $\rho\geq1.$
Let $h:\T\mapsto \R^+$ be a continuous nonnegative function satisfying \eqref{logint} and such that $h\in\cL(\T).$
Then  $O^{\rho}_{h}\in\cL $ if and only if $h$ satisfies \eqref{condition2} for some positive number $\eta\leq 1.$
\end{thm}

\begin{proof}
This proof is actually inspired from the proof of Havin-Shamoyan-Carleson-Jacobs Theorem, see for instance \cite[Theorem 3]{Dya} and \cite{HS, Shi3}.
We recall that
$$a_{h}(\theta):=\int_{|\varphi|\leq \pi\atop h(\varphi)\leq \frac{1}{2}h(\theta)}\frac{\log\frac{h(\theta)}{h(\varphi)}}{|e^{i\varphi}-e^{i\theta}|^2}d\varphi,\qquad e^{i\theta}\in\T\setminus \E_h.$$
We note that $\omega^*$ is an increasing function satisfying
\begin{eqnarray}\label{omegaprop}
\omega(\omega^*(u))=u\text{ and }\omega^*(u)\leq u,\qquad u\in[0,1].
\end{eqnarray}
For a point $e^{i\theta}\in\T\setminus \E_h$ and a point $e^{i\varphi}\in\T$ such that $|e^{i\varphi}-e^{i\theta}|\leq \omega^{*}(\frac{h(\theta)}{2\|h\|}),$
\begin{eqnarray*}
|h(\varphi)-h(\theta)|&\leq& \|h\|\omega(|e^{i\varphi}-e^{i\theta}|
\leq\frac{h(\theta)}{2},
\end{eqnarray*}
since $\omega$ is nondecreasing and satisfies \eqref{omegaprop}.
Then
\begin{eqnarray}\label{eqna}
\inf\{|e^{i\varphi}-e^{i\theta}|\ :\  h(\varphi)\leq \frac{h(\theta)}{2}\}\geq \omega^{*}(\frac{h(\theta)}{2\|h\|}).
\end{eqnarray}
In all what follows we use the notation $k_1\lesssim k_2$ to design that
for the two nonnegative functions $k_1$ and $k_2$  there exists some non specified constant $c$
such that  $k_1\leq c k_2$.
By using \eqref{eqna},
$$a_{h}(\theta)\lesssim\frac{\int_{|\varphi|\leq \pi}|\log h(\varphi)|d\varphi}{\big(\omega^{*}(\frac{h(\theta)}{2\|h\|})\big)^{2}}+\frac{|\log h(\theta)|}{\omega^{*}(\frac{h(\theta)}{2\|h\|})},$$
and hence
\begin{equation}\label{nn}
a_{h}(\theta)\lesssim \big(\omega^{*}(\frac{h(\theta)}{2\|h\|})\big)^{-2},
\end{equation}
independently of the points $e^{i\theta}\in\T\setminus \E_h$.
Now, by using \eqref{nn} and the fact that $\omega$ is nondecreasing and satisfying \eqref{cond},
\begin{eqnarray*}
 h^{\rho}(\theta)&=&(2\|h\|)^{\rho}\ \omega^{\rho}(\omega^{*}(\frac{h(\theta)}{2\|h\|}))
\\&\lesssim&(2\|h\|)^{\rho}\ \omega\Big(\big(\omega^{*}(\frac{h(\theta)}{2\|h\|})\big)^{2}\Big)
\\&\lesssim& (2\|h\|)^{\rho}\ \omega(\min\{1,a^{-1}_{h}(\theta)\}),
\end{eqnarray*}
independently of the points $e^{i\theta}\in\T\setminus \E_h.$
Therefore
\begin{equation}\label{bonjour!}
\sup_{e^{i\theta}\in\T\setminus\E_h}\frac{h^{\rho}(\theta)}{\omega(\min\{1, a^{-1}_{h}(\theta)\})}<+\infty.
\end{equation}
The desired result follows by applying Theorem \ref{asmahan1} and using the estimate \eqref{bonjour!}.
\end{proof}

We have two observations concerning Theorem \ref{coro2}. Firstly, the condition \eqref{cond} is satisfied by a large class of modulus of continuity with different growth. Especially it holds for both $\varphi_{_\beta}$ and $\psi_{_\beta}$ with $\rho=1.$
Secondly, it is known from Havin-Shamoyan-Carleson-Jacobs theorem that in the case $\omega=\chi_\alpha$  the result of Theorem \ref{coro2} is very sharp. We note that in the case $\omega=\chi_\alpha,$ we can take $\rho=2,$ and also the following estimate
\begin{equation}\label{haloo}
\sup_{\gamma\in\Gamma_{h^{2}}}\frac{|A_{h^{2}}(\gamma)|}{h^{2}(\gamma)\chi_\alpha(|\gamma|)}<+\infty,
\end{equation}
always holds, for more details on \eqref{haloo} see the proof of Theorem \ref{polycase}.

\section{ \bf An application to boundary zero sets in $\cL.$ \label{apps}}

We fix $\E\subseteq\T$ a closed subset of zero Lebesgue measure. It is clear that we can write $\T\setminus\E=\bigcup\limits_{n=0}^{\infty}\tau_{n},$ where
$\tau_n=(a_n,b_n)\subseteq\T\setminus\E$ is the open arc joining the points $a_n,b_n\in\E.$
Carleson \cite{Car} has proved that for $\E$ being a boundary zero set of a function in $\Lambda_{\chi_\alpha}$
it is necessary and sufficient that
\begin{equation*}
 \sum_{n\in\N}|\tau_n|\log |\tau_n|>-\infty.
\end{equation*}
This result was generalized  by Shirokov \cite{Shi2} for any modulus of continuity. The set $\E$ is called $\omega$-Carleson if and only if
\begin{equation}\label{rheat}
 \sum_{n\in\N}|\tau_n|\log\omega(|\tau_n|)>-\infty.
\end{equation}
Associated to $\omega$ and $\E,$ we define on $\T$ the following function
\begin{equation}\label{lfunc}
l_{\E}(e^{i\theta}):=\left\{
     \begin{array}{ll}
      \displaystyle\omega(|\tau_n|)\frac{|e^{i\theta}-b_n||e^{i\theta}-a_n|}{|\tau_n|^2}, & \qquad e^{i\theta}\in \tau_n,\\
          0,  & \qquad e^{i\theta}\in \E.\\
     \end{array}
   \right.
\end{equation}
Using the equality
$$\int_{a}^{b}\log\frac{(b-t)(t-a)}{(b-a)^2}dt=2(a-b),$$
we can check easily that $\log l_{\E}$ is integrable on $\T$  if and only if
$\E$ satisfies \eqref{rheat}. The following Theorem provide us with a simplification of some results established in \cite{Shi2}.

\begin{thm}\label{claire1}
Let $\omega$  be a modulus of continuity satisfying \eqref{cond} for some real number $\rho\geq1.$
Let $\E\subseteq\T$ be a closed subset of zero Lebesgue measure. In order that $O_{l_{\E}}^{\rho}$ belongs to $\Lambda_{\omega}$ it is necessary and sufficient
for $\E$ to be an $\omega$-Carleson set.
\end{thm}

We denote by $\mathcal{C}^{2}(O)$ the class of positive functions with continuous seconde derivatives on $O,$ an open subset of $\T.$
For proving Theorem \ref{claire1}, we need the following Lemma.

\begin{lem}\label{smouth}
Let $k:\T\mapsto \R^+$ be a continuous nonnegative function satisfying \eqref{logint}.
We suppose that $k\in\mathcal{C}^{2}(\T\setminus \E_{k})$ and that the following estimates
\begin{equation}\label{atakame1}
\sup_{e^{i\theta}\in\widetilde{\gamma}}|\frac{\partial k}{\partial \theta}(e^{i\theta})|\lesssim \frac{k(\gamma)}{\omega^{*}\big(\frac{k(\gamma)}{2\|k\|}\big)}
\qquad  \text{ and }\qquad \sup_{e^{i\theta}\in\widetilde{\gamma}}|\frac{\partial^2 k}{\partial^{2} \theta}(e^{i\theta})|\lesssim \frac{k(\gamma)}{\big(\omega^{*}\big(\frac{k(\gamma)}{2\|k\|}\big)\big)^2}
\end{equation}
are satisfied uniformly with respect to all arcs $\gamma\in\Gamma_{k}$, where $\widetilde{\gamma}$ designs the following extension of the arc $\gamma$
$$\widetilde{\gamma}:=\{e^{i(s\pm t)},\ e^{is}\in\gamma\text{ and }0\leq t\leq \omega^{*}\big(\frac{k(\gamma)}{2\|k\|}\big)\},\qquad \gamma\in\Gamma_{k}.$$
Then \eqref{condition2} holds for $h=k,$ $\eta=1$ and $\rho=1:$
\begin{equation}\label{smouthk}
 \sup_{\gamma\in\Gamma_k}\frac{|A_{k}(\gamma)|}{k(\gamma)\omega(|\gamma|)}<+\infty.
\end{equation}
\end{lem}

\begin{proof}
We let $\gamma\in\Gamma_k$ be a fixed arc.
A simple calculation gives
\begin{equation}\label{talence2}
\frac{k(\gamma)}{2}\leq k(\theta)\leq 2 k(\gamma),\qquad e^{i\theta}\in \widetilde{\gamma}.
\end{equation}
By using \eqref{atakame1} and \eqref{talence2}, and using the fact that $\omega(t)/t$ is nonincreasing
\begin{eqnarray}\nonumber
|A_{k}(\gamma)|
&=&\lim_{\varepsilon\rightarrow 0^{+}}\Big|\int_{e^{is}\in\gamma}\int_{\varepsilon}^{\omega^{*}\big(\frac{k(\gamma)}{2\|k\|}\big)} \frac{k(s+t) k(s-t)-k^{2}(s)}{t^{2}}dsdt\Big|
\\\nonumber &\leq& \int_{e^{is}\in\gamma}\int_{0}^{\omega^{*}\big(\frac{k(\gamma)}{2\|k\|}\big)}
 \Big|\frac{k(s+t) k(s-t)-k^{2}(s)}{t^{2}}\Big| dsdt
\\&\lesssim&\nonumber |\gamma|\times\omega^{*}\big(\frac{k(\gamma)}{2\|k\|}\big)\times \sup_{e^{i\theta}\in\widetilde{\gamma}}\Big(k(\gamma)\Big|\frac{\partial^{2} k}{\partial \theta^2}(e^{i\theta})\Big|
+\Big|\frac{\partial k}{\partial \theta}(e^{i\theta})\Big|^{2}\Big)
 \\&\lesssim&\label{talence3}\|k\|k(\gamma) \omega(|\gamma|).
 \end{eqnarray}
The desired estimate \eqref{smouthk} is deduced from the estimate \eqref{talence3}.
\end{proof}

\subsection*{ Proof of Theorem \ref{claire1}.}

The necessity part follows from the standard Jensen's formula. Let now prove the sufficiency part. For this, we let $\E\subseteq\T$ be an $\omega$--Carleson set for some $\omega$ satisfying \eqref{cond} with $\rho\geq1.$
We can easily check that $\|l_{\E}\|<+\infty,$ where $l_{\E}$ is the function defined by  \eqref{lfunc}. Then
$l_{\E}^{\rho}\in\cL(\T),$ and hence
$$l_{\E}^{\rho}(\gamma)\leq \|l_{\E}^{\rho}\|\omega( |\tau_n|),$$
for any arc $\gamma\subseteq\tau_n.$
Then
\begin{equation}\label{shmoz3}
\omega^{*}\big(\frac{l_{\E}^{\rho}(\gamma)}{2\|l_{\E}^\rho\|}\big)\leq|\tau_n|,
\end{equation}
for any arc $\gamma\subseteq\tau_n.$ On the other hand, we can simply calculate that
\begin{equation}\label{shmoz1}
\sup\limits_{e^{i\theta}\in\tau_{n}}\Big|\frac{\partial l_{\E}}{\partial \theta}(e^{i\theta})\Big|\lesssim \frac{\omega(|\tau_n|)}{|\tau_n|}
\qquad  \text{ and }\qquad
\sup\limits_{e^{i\theta}\in\tau_{n}}\Big|\frac{\partial^{2} l_{\E}}{\partial \theta^2}(e^{i\theta})\Big|\lesssim \frac{\omega(|\tau_n|)}{|\tau_n|^2}
\end{equation}
uniformly with respect to the numbers $n\in\N.$ Thus, by using \eqref{shmoz3} and \eqref{shmoz1} and the fact that $\omega(t)/t$ is nonincreasing,
\begin{equation*}\label{shmoz11}
\sup\limits_{e^{i\theta}\in\tau_{n}}\Big|\frac{\partial l^{\rho}_{\E}}{\partial \theta}(e^{i\theta})\Big|\lesssim \frac{l_{\E}^{\rho}(\gamma)}{\omega^{*}\big(\frac{l_{\E}^{\rho}(\gamma)}{2\|l_{\E}^\rho\|}\big)}
\qquad  \text{ and }\qquad
\sup\limits_{e^{i\theta}\in\tau_{n}}\Big|\frac{\partial^{2} l^{\rho}_{\E}}{\partial \theta^2}(e^{i\theta})\Big|\lesssim \frac{l_{\E}^{\rho}(\gamma)}{\big(\omega^{*}\big(\frac{l_{\E}^{\rho}(\gamma)}{2\|l_{\E}^\rho\|}\big)\big)^2},
\end{equation*}
uniformly with respect to the numbers $n\in\N.$
Therefore \eqref{atakame1} is satisfied for $k=l_{\E}^{\rho}.$
Hence, by applying Lemma \ref{smouth}, the condition \eqref{condition2} holds for $h=l_{\E}$ and $\eta=1.$
We now apply Theorem \ref{coro2} to deduce that $O_{l_{\E}}^{\rho}\in\Lambda_{\omega},$ which is the desired result.

\section{\bf Proof of Theorem \ref{continu}.\label{sect4}}

Let $h:\T\rightarrow \R^+$  be a nonnegative continuous function satisfying \eqref{logint}.
It is well known that the radial limits of $v_{_h}$ exist and coincide almost everywhere on $\T$ with the following function (see for instance \cite[pages 98-100]{Gar})
\begin{equation}\label{definitionv}
v_{_h}(\theta):=\lim_{\varepsilon\to 0}\frac{1}{2\pi}\int_{\varepsilon}^{\pi}\frac{\log\frac{h(\theta-t)}{h(\theta+t)}}{\tan(\frac{1}{2}t)}dt, \qquad e^{i\theta}\in\T
\setminus\E_h.
\end{equation}
The next Lemma will be used to prove the following one.

\begin{lem}\label{encorune}
Let $\theta$ and $\varphi$ be two real numbers.
We have
\begin{eqnarray}\label{regard3}
\int_{\lambda}^{\pi}\frac{\log\frac{h(\theta-t)}{h(\theta+t)}-\log\frac{h(\varphi-t)}{h(\varphi+t)}}{\tan(\frac{1}{2}t)}dt
=\frac{1}{2}\int_{\theta}^{\varphi}\int_{\lambda}^{\pi}
\frac{M_h(s,t)-M_h(s,\lambda)}{\sin^2(\frac{1}{2}t)}dsdt,
\end{eqnarray}
where $\lambda<\pi$ is a positive number and
$$M_h(s,t):=\log\frac{h(s+t)h(s-t)}{h^2(s)},$$
whenever it is well defined.
\end{lem}

\begin{proof}
A partial integration gives
\begin{eqnarray}\nonumber
&&\int_{\lambda}^{\pi}\frac{\log\frac{h(\theta-t)}
{h(\theta+t)}-\log\frac{h(\varphi-t)}{h(\varphi+t)}}{\tan(\frac{1}{2}t)}dt  \nonumber
\\\nonumber&=&\int_{\lambda}^{\pi}\frac{1}{\tan(\frac{1}{2}t)}\frac{\partial}{\partial t}\Big(\int_{\lambda}^{t}\log\frac{h(\theta-s)}{(\theta+s)}-\log\frac{h(\varphi-s)}{h(\varphi+s)}ds\Big)dt
\\\label{regard1}&=&\frac{1}{2}\int_{\lambda}^{\pi}\frac{1}{\sin^2(\frac{1}{2}t)}
\Big(\int_{\lambda}^{t}\log\frac{h(\theta-s)}{h(\theta+s)}-\log\frac{h(\varphi-s)}{h(\varphi+s)}ds\Big)dt.
\end{eqnarray}
By a change of variables
\begin{eqnarray}\nonumber
\int_{\lambda}^{t}\log\frac{h(\theta-s)}{h(\theta+s)}ds &=& -\int_{\theta-\lambda}^{\theta-t}\log h(u)du+\int_{\theta-\lambda}^{\theta-t}\log h(u+t+\lambda)du \nonumber
\\\label{suisse}&=& \int_{\theta-\lambda}^{\theta-t}
\log \frac{h(u+t+\lambda)}{h(u)}du,
\end{eqnarray}
for all numbers $t$ such that $\lambda<t\leq\pi.$
Using \eqref{suisse} and again a change of variables
\begin{eqnarray}\nonumber
&&\int_{\lambda}^{t}\Big(\log \frac{h(\theta-s)}{h(\theta+s)}-\log \frac{h(\varphi-s)}{h(\varphi+s)}\Big)ds
\\\nonumber&=& \int_{\theta-\lambda}^{\theta-t}
\log\frac{h(s+t+\lambda)}{h(s)}ds -\int_{\varphi-\lambda}^{\varphi-t}
\log\frac{h(s+t+\lambda)}{h(s)}ds
\\\nonumber&=& \int_{\theta-\lambda}^{\varphi-\lambda}
\log\frac{h(s+t+\lambda)}{h(s)}ds -\int_{\theta-t}^{\varphi-t}
\log\frac{h(s+t+\lambda)}{h(s)}ds
\\\nonumber&=&\int_{\theta}^{\varphi}\log\frac{h(s+t)}{h(s-\lambda)}ds-\int_{\theta}^{\varphi}\log\frac{h(s+\lambda)}{h(s-t)}ds
\\\label{regard2}&=&\int_{\theta}^{\varphi}\big(M_h(s,t)-M_h(s,\lambda)\big)ds.
\end{eqnarray}
By combining \eqref{regard1} and \eqref{regard2} we deduce the desired equality \eqref{regard3}.
\end{proof}

We obtain the following Lemma.

\begin{lem} \label{encorune2}
Let $\gamma:=(e^{i\theta},e^{i\varphi})\subseteq\T\setminus\E_h$ be an arc, where  $e^{i\theta},e^{i\varphi}\in\T\setminus\E_h.$
We have
\begin{eqnarray*}
&&\lim_{\varepsilon\rightarrow0}\int_{e^{is}\in\gamma}\int_{\varepsilon}^{\lambda}
\frac{M_h(s,t)}{\sin^2(\frac{1}{2}t)}dsdt
\\&=&2\lim_{\varepsilon\rightarrow0}\int_{\varepsilon}^{\lambda}\frac{\log\frac{h(\theta-t)}{h(\theta+t)}-\log\frac{h(\varphi-t)}{h(\varphi+t)}}{\tan(\frac{1}{2}t)}dt
-\int_{e^{is}\in\gamma}\int_{\lambda}^{\pi}
\frac{M_h(s,\lambda)}{\sin^2(\frac{1}{2}t)}dsdt,
\end{eqnarray*}
where $0<\lambda\leq\pi.$
\end{lem}

\begin{proof}
From Lemma \ref{encorune}
\begin{eqnarray}\nonumber\label{pain2}
&&\int_{e^{is}\in\gamma}\int_{\varepsilon}^{\lambda}
\frac{M_h(s,t)}{\sin^2(\frac{1}{2}t)}dsdt
\\\nonumber&=&2\int_{\varepsilon}^{\lambda}\frac{\log\frac{h(\theta-t)}{h(\theta+t)}-\log\frac{h(\varphi-t)}{h(\varphi+t)}}{\tan(\frac{1}{2}t)}dt
-\int_{e^{is}\in\gamma}\int_{\lambda}^{\pi}
\frac{M_h(s,\lambda)}{\sin^2(\frac{1}{2}t)}dsdt
\\&&+\int_{e^{is}\in\gamma}\int_{\varepsilon}^{\pi}
\frac{M_h(s,\varepsilon)}{\sin^2(\frac{1}{2}t)}dsdt.
\end{eqnarray}
Since $\overline{\gamma}\subseteq\T\setminus\E_h$ and $h$ is continuous
$$\lim_{\lambda\rightarrow 0} \int_{e^{is}\in\gamma}|M_h(s,\lambda)|ds=0.$$
Then, by using \eqref{regard2},
\begin{eqnarray}\nonumber
\int_{0}^{\varepsilon}\Big(\log\frac{h(\theta-s)}{h(\theta+s)}-\log\frac{h(\varphi-s)}{h(\varphi+s)}\Big)ds=\int_{e^{is}\in\gamma}M_h(s,\varepsilon)ds,
\end{eqnarray}
and hence
\begin{eqnarray}\nonumber
\lim_{\varepsilon\rightarrow 0} \frac{1}{\varepsilon}\Big|\int_{e^{is}\in\gamma}M_h(s,\varepsilon)ds\Big|
&\leq&\lim_{\varepsilon\rightarrow 0} \frac{1}{\varepsilon}\int_{0}^{\varepsilon}\Big(\Big|\log\frac{h(\theta-s)}{h(\theta+s)}\Big|+\Big|\log\frac{h(\varphi-s)}{h(\varphi+s)}\Big|\Big)ds
\\\label{pain}&=&0.
\end{eqnarray}
By combining \eqref{pain2} and \eqref{pain} we deduce the desired result.
\end{proof}

\subsection*{ \bf Proof of Theorem \ref{continu}.}

We first suppose that $O_h\in\cA(\D).$ For a fixed point $\xi\in\T\setminus\E_h,$ there exists a positive radius $r\leq1$ such that
$$\frac{1}{2}|O_h(w)|\leq|O_h(z)|\leq2|O_h(w)|,\qquad z,w\in\overline{\D(\xi,r)},$$
where
$$\D(\xi,r):=\{z\in\D\ :\ |z-\xi|\leq r\}.$$
 Then
\begin{eqnarray*}
\frac{1}{2}h(\xi)|e^{i  v_{_h}(z)}-e^{i v_{_h}(w)}|\leq |O_h(z)-O_h(w)|+|h(w)-h(z)|
\leq 2|O_h(z)-O_h(w)|,
\end{eqnarray*}
for all points $z,w\in\D(\xi,r).$
Using the following classical equality
\begin{eqnarray*}
e^{i x}-e^{i y}
&=&2\sin\big(\frac{x-y}{2}\big), \qquad x,y\in\R,
\end{eqnarray*}
we get
\begin{eqnarray}\label{continuity1}
\big|\sin\big(\frac{v_{_h}(z)-v_{_h}(w)}{2}\big)\big|
\lesssim \frac{|O_h(z)-O_h(w)|}{h(\xi)},
\end{eqnarray}
for all points $z,w\in\D(\xi,r).$ Since $O_h\in\cA(\D),$
we deduce from \eqref{continuity1} that there exists a positive radius $r'\leq r/2$ such that
\begin{eqnarray}\label{continuity2}
\big|\sin\big(\frac{v_{_h}(z)-v_{_h}(w)}{2}\big)\big|
\leq \frac{1}{2}, \qquad z,w\in\D(\xi,r').
\end{eqnarray}
Since $v_{_h}$ is continuous in $\D(\xi,r'),$ then by using \eqref{continuity2} and the standard intermediate value theorem,
\begin{eqnarray*}
|v_{_h}(z)-v_{_h}(w)|
\leq \frac{\pi}{2},
\end{eqnarray*}
for all points $z,w\in\D(\xi,r'),$ and thus
\begin{eqnarray*}
|v_{_h}(z)-v_{_h}(w)|
\lesssim \big|\sin\big(\frac{v_{_h}(z)-v_{_h}(w)}{2}\big)\big|, \qquad z,w\in\D(\xi,r').
\end{eqnarray*}
Therefore
\begin{eqnarray}\label{continuity4}
|v_{_h}(z)-v_{_h}(w)|
\lesssim \frac{|O_h(z)-O_h(w)|}{h(\xi)}, \qquad z,w\in\D(\xi,r').
\end{eqnarray}
Hence $v_{_h}$ possesses an extension to a continuous function on  $\overline{\D(\xi,r')},$
since $O_h$ is uniformly continuous on  $\overline{\D(\xi,r')}$ and by using \eqref{continuity4}.
From Lemma \ref{encorune} and Lemma \ref{encorune2}
\begin{eqnarray}\nonumber\label{continuity41}
&& v_{_h}(\theta)-v_{_h}(\varphi)
\\\nonumber&=&\lim_{\varepsilon\rightarrow 0}\int_{\varepsilon}^{\pi}\frac{\log\frac{h(\theta-t)}{h(\theta+t)}-\log\frac{h(\varphi-t)}{h(\varphi+t)}}{\tan(\frac{1}{2}t)}dt
\\\nonumber&=&\lim_{\varepsilon\rightarrow 0}
\frac{1}{2}\int_{e^{is}\gamma}\int_{\varepsilon}^{\pi}
\frac{M_h(s,t)}{\sin^2(\frac{1}{2}t)}dsdt
\\&=&\lim_{\varepsilon\rightarrow 0}
\frac{1}{2}\int_{e^{is}\in\gamma}\int_{\varepsilon}^{r'}
\frac{M_h(s,t)}{\sin^2(\frac{1}{2}t)}dsdt
+ \frac{1}{2}\int_{e^{is}\in\gamma}\int_{r'}^{\pi}
\frac{M_h(s,t)}{\sin^2(\frac{1}{2}t)}dsdt,
\end{eqnarray}
for all arcs  $\gamma:=(e^{i\theta},e^{i\varphi})\subseteq\T\setminus\E_h.$
On the other hand,
\begin{eqnarray}\label{continuity66}
\lim_{|\gamma|\rightarrow 0}\int_{e^{is}\in\gamma}\int_{r'}^{\pi}
\frac{M_h(s,t)}{\sin^2(\frac{1}{2}t)}dsdt=\lim_{|\gamma|\rightarrow 0}\int_{e^{is}\in\gamma}\int_{r'}^{\pi}
\Big|\frac{M_h(s,t)}{\sin^2(\frac{1}{2}t)}\Big|dsdt=0,
\end{eqnarray}
uniformly with respect to all arcs  $\gamma\subseteq\overline{\D(\xi,r')}\cap\T,$
since
\begin{eqnarray*}
\int_{e^{is}\in\gamma}\int_{r'}^{\pi}
\Big|\frac{M_h(s,t)}{\sin^2(\frac{1}{2}t)}\Big|dsdt\lesssim|\gamma|\Big(
\frac{\int_{|\varphi|\leq \pi}|\log h(\varphi)|d\varphi}{(r^{'})^{2}}+\frac{|\log2 h(\xi)|}{r'}\Big).
\end{eqnarray*}
By using \eqref{continuity41}, \eqref{continuity66} and the fact that $v_{_h}$ is uniformly continuous on $\overline{\D(\xi,r')}\cap\T,$ we obtain
\begin{eqnarray}\label{continuity5}
\lim_{|\gamma|\rightarrow 0}\lim_{\varepsilon\rightarrow 0}
\int_{e^{is}\in\gamma}\int_{\varepsilon}^{r'}
\frac{M_h(s,t)}{\sin^2(\frac{1}{2}t)}dsdt=0,
\end{eqnarray}
uniformly with respect to all arcs  $\gamma\subseteq\overline{\D(\xi,r')}\cap\T.$
We also note that
\begin{equation}\label{continuity6}
\frac{M_h(s,t)}{\sin^2(\frac{1}{2}t)}\asymp \frac{h(s+t) h(s-t)-h^2(s)}{t^2h^{2}(\xi) },\qquad e^{is}\in\overline{\D(\xi,r')}\text{ and }t\leq r'.
\end{equation}
Hence \eqref{continuity} is deduced from \eqref{continuity5} and \eqref{continuity6} with $\lambda=r'$ and $\mathcal{U}=\overline{\D(\xi,r')}\cap\T.$

Let now prove the sufficiency part. We note first that $O_h\in\cH^\infty(\D)$ and that $|O_h|$ is continuous on $\overline{\D},$ since $h$ is continuous.
Now, since a harmonic extension over $\D$ of a continuous function on $\T$ is continuous on $\overline{\D},$ to prove that $O_h\in\cA(\D)$ we just need to show that if \eqref{continuity} holds then the radial limits of $O_h$ exist on $\T\setminus\E_h,$ and this limits define a continuous function on $\T$ with zero values on $\E_h.$ To prove that it is sufficient to show that the limits $v_{_h},$
defined in \eqref{definitionv}, are finite and define a continuous function on $\T\setminus\E_h.$
We suppose that  \eqref{continuity} is satisfied on a neighborhood $\mathcal{U}$ of a point $\xi\in\T\setminus\E_h,$ where  $0<\lambda\leq\pi.$
Since $h$ is continuous, there exists a positive radius $r\leq\pi$ such that $\overline{\D(\xi,r)}\cap\T\subseteq\mathcal{U}$ and
$$\frac{1}{2}h(\theta)\leq h(\varphi)\leq 2 h(\theta),\qquad \text{ for all }e^{i\theta},e^{i\varphi}\in\overline{\D(\xi,r)}.$$
Then, as we have done for computing \eqref{continuity6},
\begin{equation*}
\frac{M_h(s,t)}{\sin^2(\frac{1}{2}t)}\asymp \frac{h(s+t) h(s-t)-h^2(s)}{t^2h^{2}(\xi) },\qquad e^{is}\in\overline{\D(\xi,\lambda')}\text{ and }t\leq \lambda',
\end{equation*}
where  $\lambda':=\min\{\lambda,r/2\}.$
Thus, by using the hypothesis \eqref{continuity},
\begin{eqnarray}\label{continuity7}
\lim_{|\gamma|\rightarrow 0}\lim_{\varepsilon\rightarrow 0}
\int_{e^{is}\in\gamma}\int_{\varepsilon}^{\lambda'}
\frac{M_h(s,t)}{\sin^2(\frac{1}{2}t)}dsdt=0
\end{eqnarray}
uniformly with respect to all arcs $\gamma\subseteq\overline{\D(\xi,\lambda')}\cap\T.$
On the other hand, as we have calculated in \eqref{continuity66},
\begin{eqnarray}\label{continuity8}
\lim_{|\gamma|\rightarrow 0}\int_{e^{is}\in\gamma}\int_{\lambda'}^{\pi}
\frac{M_h(s,t)}{\sin^2(\frac{1}{2}t)}dsdt=0,
\end{eqnarray}
uniformly with respect to all arcs  $\gamma\subseteq\overline{\D(\xi,\lambda')}\cap\T.$
From \eqref{continuity41}, \eqref{continuity7} and \eqref{continuity8} and the well known fact that the set of all points $e^{i\varphi}\in\T$ for which the limit $v_{_h}(\varphi)$ is finite is dense on $\T,$ we deduce that $v_{_h}$ possesses an extension to a continuous function on $\overline{\D(\xi,\lambda')}\cap\T.$  Therefore $v_{_h}$ is continuous on $\T\setminus\E_h,$ and hence $O_h$ is continuous on $\T.$  This completes the proof of Theorem \ref{continu}.

\section{\bf Proof of Theorem \ref{polycase} \label{proofpoly}.}

By applying Lemma \ref{linda1} and Lemma \ref{linda2} below we deduce that if $O_h\in\cL$ then  \eqref{second} holds, and if \eqref{second} holds then so is for \eqref{first}.
We now suppose that \eqref{first} is satisfied.
By applying Theorem \ref{asmahan1}, to prove that $O_h\in\cL$ it is sufficient to show that \eqref{condition2} holds for $\rho=1$ and $\eta=1$ provided that \eqref{c} is verified.
We let $\widetilde{\gamma}$ be the following extension of a fixed arc $\gamma:=(e^{i\theta},e^{i\varphi})\in\Gamma_{h}$
\begin{equation}\label{bienne00}
\widetilde{\gamma}:=\{e^{i(s\pm t)},\ e^{is}\in\gamma\text{ and }0\leq t\leq \omega^{*}\big(\frac{h(\gamma)}{2\|h\|}\big)\}.
\end{equation}
Since
\begin{equation}\label{bienne0}
\frac{h(\gamma)}{2}\leq h(\theta)\leq 2 h(\gamma),\qquad e^{i\theta}\in \widetilde{\gamma},
\end{equation}
then, for every point $e^{is}\in\gamma$  and  every number $t$ such that $0< t\leq \omega^{*}\big(\frac{h(\gamma)}{2\|h\|}\big),$
\begin{equation}\label{bienne}
\frac{h(s+t) h(s-t)-h^2(s)}{t^2h^{2}(\gamma)}\asymp  \frac{M_h(s,t)}{\sin^2(\frac{1}{2}t)}.
\end{equation}
Thus
\begin{eqnarray}\label{last1}\nonumber
|A_{h}(\gamma)|
&\lesssim& \Big|\lim_{\varepsilon\rightarrow 0}\int_{e^{is}\in\gamma}\int_{\varepsilon}^{|\gamma|}
 \frac{h(s+t) h(s-t)-h^2(s)}{t^{2}}dsdt\Big|
\\\nonumber&&+\int_{e^{is}\in\gamma}\int_{|\gamma|}^{\omega^{*}(\frac{h(\gamma)}{2\|h\|})}
 \Big|\frac{h(s+t) h(s-t)-h^2(s)}{t^{2}}\Big|dsdt
 \\\nonumber&\lesssim& h^{2}(\gamma)\Big|\lim_{\varepsilon\rightarrow 0}\int_{e^{is}\in\gamma}\int_{\varepsilon}^{|\gamma|}
\frac{M_h(s,t)}{\sin^2(\frac{1}{2}t)}dsdt\Big|
\\&&+\|h\|h(\gamma)|\gamma|\int_{|\gamma|}^{1}
 \frac{\omega(t)}{t^{2}}dt.
\end{eqnarray}
From Lemma \ref{encorune2}
\begin{eqnarray*}
&&\Big|\lim_{\varepsilon\rightarrow0}\int_{e^{is}\in\gamma}\int_{\varepsilon}^{|\gamma|}
\frac{M_h(s,t)}{\sin^2(\frac{1}{2}t)}dsdt\Big|
\\&=&\Big|2\lim_{\varepsilon\rightarrow0}\int_{\varepsilon}^{|\gamma|}\frac{\log\frac{h(\theta-t)}{h(\theta+t)}-\log\frac{h(\varphi-t)}{h(\varphi+t)}}{\tan(\frac{1}{2}t)}dt
-\int_{e^{is}\in\gamma}\int_{|\gamma|}^{\pi}
\frac{M_h(s,|\gamma|)}{\sin^2(\frac{1}{2}t)}dsdt\Big|
\\&\lesssim&\int_{0}^{|\gamma|}\frac{|\log\frac{h(\theta-t)}{h(\theta+t)}|+|\log\frac{h(\varphi-t)}{h(\varphi+t)}|}{t}dt
+\int_{e^{is}\in\gamma}\int_{|\gamma|}^{\pi}
\frac{|M_h(s,|\gamma|)|}{t^2}dsdt.
\end{eqnarray*}
A simple calculation gives(with \eqref{bienne0} in mind)
\begin{eqnarray*}\label{last4}\nonumber
|\log\frac{h(s-t)}{h(s+t)}|\lesssim \frac{|h(s-t)-h(s+t)|}{h(\gamma)}
\lesssim \frac{\|h\|\omega(t)}{h(\gamma)},
\end{eqnarray*}
for all points $e^{is}\in\overline{\gamma}$  and all numbers $t$ such that $0<t\leq\omega^{*}\big(\frac{h(\gamma)}{2\|h\|}\big).$
Also
$$M_h(s,|\gamma|)\lesssim \frac{\|h\|\omega(|\gamma|)}{h(\gamma)},\qquad e^{is}\in\gamma.$$
Therefore
\begin{eqnarray}\label{last2}
\Big|\lim_{\varepsilon\rightarrow 0}\int_{e^{is}\in\gamma}\int_{\varepsilon}^{|\gamma|}
\frac{M_h(s,t)}{\sin^2(\frac{1}{2}t)}dsdt\Big|\lesssim\frac{\|h\|}{h(\gamma)} \int_{0}^{|\gamma|}\frac{\omega(t)}{t}dt+ \frac{\|h\|\omega(|\gamma|)}{h(\gamma)}.
\end{eqnarray}
Since $\omega$ verifies \eqref{c}, we deduce from \eqref{last1} and \eqref{last2} the following estimate
\begin{eqnarray}\label{last3}\nonumber
|A_{h}(\gamma)|&\lesssim& \|h\|h(\gamma)\Big(\int_{0}^{|\gamma|}\frac{\omega(t)}{t}dt + |\gamma|\int_{|\gamma|}^{1}
 \frac{\omega(t)}{t^{2}}dt+\omega(|\gamma|)\Big)
 \\&\lesssim& \|h\|h(\gamma)\omega(|\gamma|).
\end{eqnarray}
Thus \eqref{condition2} holds, and hence the proof of Theorem \ref{polycase} is completed.

\section{\bf Proof of Theorem \ref{asmahan1}.\label{proftheorem} }

Let $h:\T\rightarrow \R^+$  be a nonnegative continuous function satisfying \eqref{logint} and  such that  $h\in\cL(\T)$ for some modulus of continuity $\omega.$ The proof of Theorem \ref{asmahan1} will be deduced from a series of technical lemmas.

\subsection{\bf The Necessity part.\label{ver}}

We begin by the following lemma, which gives rise to the necessary condition \eqref{condition1} of Theorem \ref{asmahan1}.

\begin{lem}\label{linda1} Let $\rho\geq1$ be a number.
We suppose that
\begin{eqnarray}\label{linda00}
 \sup_{z\in\D\atop 1-|z|=\omega^{*}(\frac{\delta h^{\rho}(z/|z|)}{2\|h\|^{\rho}})}\Big|\log\frac{|O_{h}(z)|}{h(z/|z|)}\Big|<+\infty,
\end{eqnarray}
for some positive number $\delta\leq1.$ Then \eqref{condition1} is satisfied.
\end{lem}

\begin{proof}  We recall that
$$\log|O_{h}(z)|=u_{_h}(z)=\frac{1}{2\pi}\int_{-\pi}^{\pi}\frac{1-|z|^2}{|e^{i\varphi}-z|^2}\log h(\varphi)d\varphi, \qquad z\in\D.$$
We let $e^{i\theta}\in\T\setminus \E_h$ be a fixed point such that $a^{-1}_{h}(\theta)\leq 1.$  We have
\begin{eqnarray}\nonumber
&&\frac{1}{2\pi}\int_{[-\pi,\pi[\setminus\Omega_{h}(\theta)}\frac{1-|z|^2}{|e^{i\varphi}-z|^2}\log\frac{h(\varphi)}{h(\theta)}d\varphi
\\\label{linda11}&=&\log\frac{|O_{h}(z)|}{h(\theta)}-\frac{1}{2\pi}\int_{\Omega_{h}(\theta)}\frac{1-|z|^2}{|e^{i\varphi}-z|^2}\log\frac{h(\varphi)}{h(\theta)}d\varphi,
\qquad z\in\D,
\end{eqnarray}
where $$\Omega_{h}(\theta):=\{\varphi\in[-\pi,\pi[\ :\ |e^{i\varphi}-e^{i\theta}|\leq\omega^{*}\big(\frac{h(\theta)}{2\|h\|}\big):=\lambda(\theta)\},\qquad e^{i\theta}\in\T.$$
We also have
\begin{eqnarray} \label{linda111}
\frac{1}{2\pi}\int_{\Omega_{h}(\theta)}\frac{1-|z|^2}{|e^{i\varphi}-z|^2}\big|\log\frac{h(\varphi)}{h(\theta)}\big|d\varphi
\leq  \frac{\log 2}{2\pi}\int_{\Omega_{h}(\theta)}\frac{1-|z|^2}{|e^{i\varphi}-z|^2}d\varphi\leq \log 2, \quad z\in\D.
\end{eqnarray}
For the following point
$$z(\theta):=\big(1-\omega^{*}(\frac{\delta h^{\rho}(\theta)}{2\|h\|^{\rho}})\big)e^{i\theta},$$
we simply compute that
\begin{eqnarray} \label{linda1010}
0<1-|z(\theta)|=\omega^{*}(\frac{\delta h^{\rho}(\theta)}{2\|h\|^{\rho}})\leq \lambda(\theta)\leq\frac{1}{2}.
\end{eqnarray}
By using \eqref{linda00}, \eqref{linda11} and \eqref{linda111}
\begin{eqnarray}\label{linda12}
\Big|\frac{1}{2\pi}\int_{[-\pi,\pi[\setminus\Omega_{h}(\theta)}\frac{1-|z(\theta)|^2}{|e^{i\varphi}-z(\theta)|^2}\log\frac{h(\varphi)}{h(\theta)}d\varphi\Big|
\leq c,
\end{eqnarray}
where $c>0$ is a constant not depending on the points $e^{i\theta}\in\T\setminus \E_h.$
We use the following classical equality
\begin{equation}\label{euclid}
|z-w|^2=||z|-|w||^2+|zw||\frac{z}{|z|}-\frac{w}{|w|}|^2,\qquad z,w\in\overline{\D},
\end{equation}
to calculate
\begin{eqnarray}\label{linda13}
\frac{1}{2}|e^{i\varphi}-e^{i\theta}|^2\leq|e^{i\varphi}-z(\theta)|^2\leq 2|e^{i\varphi}-e^{i\theta}|^2,\qquad \varphi\in[-\pi,\pi[\setminus\Omega_{h}(\theta).
\end{eqnarray}
Now, using the fact that $\omega(t)/t$ is nonincreasing,
\begin{eqnarray}\nonumber
\big|\log \frac{h(\varphi)}{h(\theta)}\big|&\leq&\log\big(2\|h\|\frac{\omega(|e^{i\varphi}-e^{i\theta}|)}{h(\theta)}+1\big)
\\\nonumber&\leq&\log\big(\frac{|e^{i\varphi}-e^{i\theta}|}{\lambda(\theta)}+1\big),
\end{eqnarray}
for all points $\varphi\in[-\pi,\pi[\setminus\Omega_{h}(\theta)$ such that  $h(\varphi)\geq\frac{1}{2} h(\theta).$
Thus, with \eqref{linda1010} and \eqref{linda13} in mind,
\begin{eqnarray}\nonumber
&&\int_{\varphi\in[-\pi,\pi[\setminus\Omega_{h}(\theta)\atop h(\varphi)\geq\frac{1}{2} h(\theta)}\frac{1-|z(\theta)|^2}{|e^{i\varphi}-z_h(\theta)|^2}\big|\log\frac{h(\varphi)}{h(\theta)}\big|d\varphi
\\\nonumber&\leq&4\int_{|e^{i\varphi}-e^{i\theta}|\geq\lambda(\theta)}\frac{\lambda(\theta)}{|e^{i\varphi}-e^{i\theta}|^2}\log\big(\frac{|e^{i\varphi}-e^{i\theta}|}
{\lambda(\theta)}+1\big)d\varphi
\\\label{linda14}&\lesssim& \int_{t\geq1}\frac{\log t}{t^2}dt.
\end{eqnarray}
Therefore
\begin{eqnarray}\nonumber\label{linda15}
\omega^{*}(\frac{\delta h^{\rho}(\theta)}{2\|h\|^{\rho}}) a_{h}(\theta)
&=&\int_{\varphi\in[-\pi,\pi[\atop h(\varphi)\leq\frac{1}{2} h(\theta)}\frac{1-|z(\theta)|}{|e^{i\varphi}-e^{i\theta}|^2} \log\frac{h(\theta)}{h(\varphi)}d\varphi
\\\nonumber&\leq&2\int_{\varphi\in[-\pi,\pi[\atop h(\varphi)\leq\frac{1}{2} h(\theta)}\frac{1-|z(\theta)|^2}{|e^{i\varphi}-z(\theta)|^2}\log\frac{h(\theta)}{h(\varphi)}d\varphi
\\\nonumber&\leq& 2\Big|\int_{[-\pi,\pi[\setminus\Omega_{h}(\theta)}\frac{1-|z(\theta)|^2}{|e^{i\varphi}-z(\theta)|^2}\log\frac{h(\varphi)}{h(\theta)}d\varphi\Big|
\\&&+ 2\int_{\varphi\in[-\pi,\pi[\setminus\Omega_{h}(\theta)\atop h(\varphi)\geq\frac{1}{2} h(\theta)}\frac{1-|z(\theta)|^2}{|e^{i\varphi}-z(\theta)|^2}\big|\log\frac{h(\varphi)}{h(\theta)}\big|d\varphi.
\end{eqnarray}
From \eqref{linda12}, \eqref{linda14} and \eqref{linda15} we deduce the desired result.
\end{proof}

The necessity of the condition \eqref{condition1} so that $O^{\rho}_h\in\cL$ is deduced by combining Lemma \ref{linda1} with the following lemma.

\begin{lem}\label{linda2}
We suppose that $O^{\rho}_{h}\in\cL,$ where $\rho\geq1.$  Then there exists a positive number $\delta\leq 1$ satisfying \eqref{linda00}.
\end{lem}

\begin{proof}
We suppose that $O^{\rho}_{h}\in\cL.$  For a point $z\in\D$ such that $1-|z|=\omega^{*}(\frac{\delta h^{\rho}(z/|z|)}{2\|h\|^{\rho}}),$
where $\delta:=\min\{1,\frac{\|h\|^{\rho}}{\|O^{\rho}_h\|}\},$ we simply compute that
$$\big||O^{\rho}_h(z)|-h^{\rho}(z/|z|)\big|\leq \|O^{\rho}_h\|\omega(1-|z|)\leq\frac{h^{\rho}(z/|z|)}{2},$$
and hence
$$\frac{1}{2}h^{\rho}(z/|z|)\leq |O^{\rho}_h(z)|\leq\frac{3}{2}h^{\rho}(z/|z|).$$
It follows
\begin{eqnarray*}
 \sup_{z\in\D \atop 1-|z|=\omega^{*}(\frac{\delta h^{\rho}(z/|z|)}{2\|h\|^{\rho}})}\Big|\log\frac{|O_{h}(z)|}{h(z/|z|)}\Big|\leq\frac{\log2}{\rho},
\end{eqnarray*}
which is the desired result.
\end{proof}

The next Lemma shows that \eqref{condition2} is a necessary condition so that $O^{\rho}_h\in\cL.$

\begin{lem}\label{linda3}
We suppose that $O_{h}\in\cL.$ Then there exists a positive number $\eta\leq 1$ such that
\begin{equation*}
 \sup_{\gamma\in\Gamma_{h,\eta}}\frac{|A_h(\gamma)|}{h(\gamma)\omega(|\gamma|)}<+\infty.
 \end{equation*}
\end{lem}

\begin{proof}
We fix two distinct points $e^{i\theta}$ and $e^{i\varphi}$ such that $-\pi\leq\theta<\varphi<\pi.$ We suppose that
$\gamma\in\Gamma_h,$ where $\gamma=(e^{i\theta},e^{i\varphi}):=\{e^{is}\ :\ \theta< s<\varphi\}.$
We set
\begin{equation}\label{biel}
\mu:=\inf\limits_{e^{is}\in\gamma}a^{-1}_{h}(s).
\end{equation}
Since  $\gamma\in\Gamma_h,$ then $\lambda:=\omega^{*}(\frac{h(\gamma)}{2\|h\|})>0,$ and hence $h(\gamma)>0.$ It follows that
$\overline{\gamma}\subseteq\T\setminus \E_h$ and by consequence $\mu>0.$
We have the following two possible cases.

$\mathbf{A}.$ We first suppose that $|\gamma|\geq \mu.$ Since $\gamma\in\Gamma_h,$ then $\mu\leq\lambda\leq \omega^{*}(\frac{1}{2})\leq \frac{1}{2}.$
Since $O_{h}\in\cL$ then by applying
both Lemma \ref{linda1} and Lemma \ref{linda2} we deduce that \eqref{condition1} holds, and hence
\begin{equation}\label{bordeaux1}
 h(\gamma)\leq h(s_{_0})\lesssim \omega(a^{-1}_{h}(s_{_0}))
=  \omega(\mu) \leq   \omega(|\gamma|),
\end{equation}
where $s_{_0}$ is a real number such that $e^{is_{_0}}\in\overline{\gamma}$ and
$\mu=a^{-1}_{h}(s_{_0}).$

$\mathbf{B}.$ We now suppose that $|\gamma|\leq\min\{\lambda,\ \mu\}.$ Since $\overline{\gamma}\subseteq\T\setminus \E_h$  and $O_{h}\in\cL,$ then  the limits $v_{_h}(\theta)$ and $v_{_h}(\varphi)$ are finite and
\begin{eqnarray}\label{aurelie1}
h(\gamma)|e^{i  v_{_h}(\theta)}-e^{i v_{_h}(\varphi)}|
\lesssim |O_h(\theta)-O_h(\varphi)|+|h(\varphi)-h(\theta)|\lesssim
 \|O_h\| \omega(|\gamma|).
\end{eqnarray}
Moreover, by using  Lemma \ref{encorune2},
\begin{eqnarray}\nonumber
e^{i v_{_h}(\theta)}-e^{i v_{_h}(\varphi)}
&=&2\sin\Big(\frac{v_{_h}(\theta)-v_{_h}(\varphi)}{2}\Big)
\\\nonumber&=&2\sin(I_{h}(\gamma)+J_{h}(\gamma))
\\\label{aurelie1.5}&=&2\sin(I_{h}(\gamma))
\cos(J_{h}(\gamma))
+2\cos(I_{h}(\gamma))\sin(J_{h}(\gamma)),
\end{eqnarray}
where
$$I_{h}(\gamma):=\lim_{\varepsilon\rightarrow 0}\frac{1}{8\pi}\int_{\theta}^{\varphi}\int_{\varepsilon}^{\lambda}
\frac{M_h(s,t)}{\sin^2(\frac{1}{2}t)}dsdt$$
and
$$J_{h}(\gamma):=\frac{1}{8\pi}\int_{\theta}^{\varphi}\int_{\lambda}^{\pi}
\frac{M_h(s,t)}{\sin^2(\frac{1}{2}t)}dsdt.$$
Then, by using the following simple inequality,
 $$|x|+|\cos x|\geq 1,\qquad x>0,$$
we obtain
\begin{eqnarray}\nonumber
\big|\sin(I_{h}(\gamma))\big|
&\leq&\big|\sin(I_{h}(\gamma))J_{h}(\gamma)\big|
+\big|\sin(I_{h}(\gamma))\cos(J_{h}(\gamma))\big|
\\\nonumber&\leq&\big|J_{h}(\gamma)\big|
+ \frac{1}{2}\big|e^{i v_{_h}(\theta)}-e^{i v_{_h}(\varphi)}\big|
+\big|\cos(I_{h}(\gamma))\sin(J_{h}(\gamma))\big|
\\\label{aurelie2}&\leq& \frac{1}{2} \big|e^{i v_{_h}(\theta)}-e^{i v_{_h}(\varphi)}\big| + 2\big|J_{h}(\gamma)\big|.
\end{eqnarray}
From \eqref{aurelie1} and \eqref{aurelie2}
\begin{eqnarray}\label{bateau1}
h(\gamma)\big|\sin(I_{h}(\gamma))\big|
\lesssim\|O_h\|\omega(|\gamma|)+ h(\gamma)\big|J_{h}(\gamma)\big|.
\end{eqnarray}
We have
\begin{eqnarray}\nonumber
&&h(\gamma)|J_{h}(\gamma)|
\\\nonumber &\lesssim& h(\gamma)\int_{\theta}^{\varphi}\int_{\lambda\leq|t|\leq\pi}
\frac{|\log\frac{h(s+ t)}{h(s)}|}{\sin^2(\frac{1}{2}t)}dsdt
\\\nonumber&=&h(\gamma)\int_{\theta}^{\varphi}\Big(\int_{\lambda\leq|t|\leq\pi\atop h(s+ t)\leq \frac{1}{2}h(s)}
\frac{\log\frac{h(s)}{h(s+t)}}{\sin^2(\frac{1}{2}t)}dt \Big)ds
\\\nonumber&&+  h(\gamma)\int_{\theta}^{\varphi}\Big(\int_{\lambda\leq|t|\leq\pi\atop h(s+ t)\geq \frac{1}{2}h(s)}
\frac{|\log\frac{h(s+ t)}{h(s)}|}{\sin^2(\frac{1}{2}t)}dt \Big)ds
\\\label{bateau3}&\leq& \int_{\theta}^{\varphi} h(s)a_{h}(s)ds
+ h(\gamma)\int_{\theta}^{\varphi}\Big(\int_{\lambda\leq|t|\leq\pi\atop h(s+ t)\geq \frac{1}{2}h(s)}
\frac{|\log\frac{h(s+ t)}{h(s)}|}{t^2}dt \Big)ds.
\end{eqnarray}
For all points $e^{is}\in\gamma$ and all numbers $t$ such that $|t|\geq\lambda$ and such that $h(s+t)\geq \frac{1}{2}h(s),$ we have
\begin{eqnarray}\nonumber
\big|\log\frac{h(s+t)}{h(s)}\big|&=&\log\big(\frac{\big|h(s+ t)-h(s)\big|}{\inf\{h(s+ t),h(s)\}}+1\big)
\\\nonumber&\leq&\log\big(\frac{2\|h\|}{ h(\gamma)}\omega(|e^{i (s+ t)}-e^{is}|)+1\big)
\\\nonumber&\leq& \log\big(\frac{2\|h\|\omega(|t|)}{ h(\gamma)}+1\big),
\end{eqnarray}
then, by using the fact that  $\omega(t)/t$ is nonincreasing and \eqref{omegaprop},
\begin{eqnarray}\label{3yit}
\big|\log\frac{h(s+ t)}{h(s)}\big|\leq\log\big(\frac{2\|h\|\omega(\lambda)}{h(\gamma)}\frac{|t|}{\lambda}+1\big)
=\log\big(\frac{|t|}{\lambda}+1\big).
\end{eqnarray}
Thus, by using \eqref{3yit} and again the fact that $\omega(t)/t$ is nonincreasing,
\begin{eqnarray}\nonumber
h(\gamma)\int_{\theta}^{\varphi}\Big(\int_{\lambda\leq|t|\leq\pi\atop h(s+ t)\geq \frac{1}{2}h(s)}
\frac{|\log\frac{h(s+ t)}{h(s)}|}{t^2}dt \Big)ds
&\lesssim&  h(\gamma)|\gamma|\int_{t\geq \lambda}
\frac{\log\big(\frac{t}{\lambda}+1\big)}{t^2}dt
\\\nonumber&\lesssim&  \frac{|\gamma|h(\gamma)}{\lambda}\int_{u\geq1}\frac{\log (u+1)}{u^2}du
\\\label{3yit2}&\lesssim& \|h\|\omega(|\gamma|).
\end{eqnarray}
By using the estimate \eqref{condition1} and
also once again the fact that $\omega(t)/t$ is nonincreasing,
\begin{eqnarray}\label{etencore}
\int_{\theta}^{\varphi}h(s)a_{h}(s)ds
\lesssim  |\gamma|\sup\limits_{e^{is}\in\gamma} h(s)a_{h}(s)
\lesssim |\gamma|\frac{\omega(\min\{1,\mu\})}{\min\{1,\mu\}}
\lesssim  \omega(|\gamma|).
\end{eqnarray}
From \eqref{bateau3}, \eqref{3yit2} and \eqref{etencore} we deduce
\begin{eqnarray}\label{bateau2}
h(\gamma)|J_{h}(\gamma)| \lesssim\omega(|\gamma|).
\end{eqnarray}
Therefore, from \eqref{bateau1} and \eqref{bateau2},
\begin{equation}\label{bordeaux2}
h(\gamma)\big|\sin(I_{h}(\gamma))\big|
\lesssim\omega(|\gamma|).
\end{equation}
Hence, by combining  the above two situations $\mathbf{A}$ and $\mathbf{B},$
\begin{equation}\label{bordeaux3}
h(\gamma)\big|\sin\big(I_{h}(\gamma)\big)\big|\leq M_h\omega(|\gamma|),\qquad \gamma\in\Gamma_h.
\end{equation}
where $M_h>0$ is a constant not depending of the arcs $\gamma\in\Gamma_h.$

We will now deduce from \eqref{bordeaux3} and the fact that $O_h\in\cA(\D)$ an estimation of
$I_{h}(\gamma)$ without the function $\sin.$
Let  $\eta:=\min\{1,\ \|h\|/M_h\}.$ From \eqref{bordeaux3} we get
\begin{equation}\label{bordeaux4}
\big|\sin\big(I_{h}(\gamma)\big)\big|\leq \frac{1}{2}, \qquad  \gamma\in\Gamma_{h,\eta}.
\end{equation}
As above in \eqref{bienne},
\begin{equation} \label{biel0}
\frac{M_h(s,t)}{\sin^2(\frac{1}{2}t)}\asymp \frac{h(s+t)h(s-t)-h^{2}(s)}{t^{2}h^{2}(\gamma)},\qquad \theta\leq s \leq \varphi\text{ and } 0<t\leq\lambda.
\end{equation}
Then
\begin{equation}\label{bordeaux66}
I_{h}(\gamma)\asymp \frac{A_{h}(\gamma)}{h^{2}(\gamma)},\qquad \gamma\in\Gamma_{h}.
\end{equation}
We fix an arc $\gamma_0\in\Gamma_{h,\eta}.$ We observe $\overline{\gamma_0}\subseteq\T\setminus\E_h$ and that for all arcs $\gamma$ such that $\gamma\subseteq\gamma_0$ we have $\gamma\in\Gamma_{h,\eta}.$  Since $O_h\in\cA(\D)$ then $h$ satisfies \eqref{continuity} and hence
$A_{h}(\gamma)$ goes to 0 when $|\gamma|$ tends to 0 uniformly
with respect to all arcs $\gamma\subseteq\gamma_0.$ Thus  $I_{h}(\gamma)$ also goes to 0 when $|\gamma|$ tends to 0 uniformly
with respect to all arcs $\gamma\subseteq\gamma_0,$ since it satisfies \eqref{bordeaux66}.
By applying the classical intermediate value theorem  and using \eqref{bordeaux4}, we obtain
$$\big|I_{h}(\gamma)\big|\leq \frac{\pi}{2},\qquad \gamma\in\Gamma_{h,\eta}.$$
It follows
\begin{equation}\label{bordeaux5}
\big|I_{h}(\gamma)\big|\lesssim \big|\sin\big(I_{h}(\gamma)\big)\big|, \qquad  \gamma\in\Gamma_{h,\eta}.
\end{equation}
Therefore, by \eqref{bordeaux3} and \eqref{bordeaux5},
\begin{equation*}\label{bordeaux6}
h(\gamma)\big|I_{h}(\gamma)\big|\lesssim M_h\ \omega(|\gamma|),\qquad  \gamma\in\Gamma_{h,\eta}.
\end{equation*}
Hence
\begin{equation}\label{bordeaux6}
\frac{\big|A_{h}(\gamma)\big|}{h(\gamma)}\lesssim M_h \omega(|\gamma|),\qquad \gamma\in\Gamma_{h,\eta}.
\end{equation}
Which finishes the proof of Lemma \ref{linda3}.
\end{proof}

\subsection{\bf The sufficiency part.}

Our last Lemma is the following.

\begin{lem}\label{construire2}
If the conditions \eqref{condition1} and \eqref{condition2} hold, then   $O^{\rho}_h\in\cL.$
\end{lem}

\begin{proof}
We note that $O^{\rho}_h=O_{h^{^\rho}}.$
Since $h^\rho$ satisfies \eqref{condition2} then we can simply check that also it satisfies   \eqref{continuity}, by considering the above estimates \eqref{continuity66} and \eqref{continuity6}.
Thus $O^{\rho}_h\in\cA(\D),$ by applying Theorem \ref{continu}. According to Tamrazov's theorem \cite{Tam}, to prove that $O^{\rho}_h\in\cL$ it is sufficient to show that
$O^{\rho}_h\in\cL(\T),$ see also \cite[Appendix A]{Bou}.  Let  $e^{i\theta}\in\T$ and $e^{i\varphi}\in\T$
be two distinct points such that $-\pi\leq\theta<\varphi<\pi.$  We have
\begin{eqnarray}\nonumber
|O^{\rho}_h(\theta)-O^{\rho}_h(\varphi)|
&\lesssim&  \|h^{\rho}\|\ \omega(|\gamma|)+ h^{\rho}(\gamma)|e^{i\rho v_{_h}(\theta)}-e^{i\rho v_{_h}(\varphi)}|
\\\label{sift}&\lesssim&   \omega(|\gamma|)+ \min\{2h^{\rho}(\gamma),\ h^{\rho}(\gamma)|e^{i\rho v_{_h}(\theta)}-e^{i\rho v_{_h}(\varphi)}|\},
\end{eqnarray}
where $\gamma:=(e^{i\theta},e^{i\varphi}).$
We have the following situations.

$\mathbf{A}.$ We suppose that  $|\gamma|\geq \lambda,$ where $\lambda:=\omega^{*}(\frac{\eta h^\rho(\gamma)}{2\|h^\rho\|}).$
Then
\begin{eqnarray*}
h^{\rho}(\gamma)=\frac{2}{\eta}\|h^{\rho}\|\omega(\lambda)\leq \frac{2}{\eta}\|h^{\rho}\| \omega(|\gamma|).
\end{eqnarray*}

$\mathbf{B}.$ We now suppose that $|\gamma|\leq\lambda.$ Then $\gamma\in\Gamma_{h^{\rho},\eta},$ $\lambda>0$ and hence $h(\gamma)>0.$
It follows $\overline{\gamma}\subseteq\T\setminus \E_h$ and thus the number
 \begin{equation}
\mu:=\inf\limits_{e^{is}\in\gamma}a^{-1}_{h}(s)
\end{equation}
is positive.
We let $e^{is_{_0}}\in\overline{\gamma}$ be a point such that
$a^{-1}_{h}(s_{_0})=\mu.$ We consider the following two cases.
\begin{itemize}
  \item [$B_1.$] If $|\gamma|\geq \mu,$ then $\mu\leq\lambda\leq\frac{1}{2}.$ By using the hypothesis \eqref{condition1}
\begin{eqnarray*}
h^{\rho}(\gamma)\leq h^{\rho}(s_{_0})\lesssim  \omega(a^{-1}_{h}(s_{_0}))=  \omega(\mu)\leq   \omega(|\gamma|).
\end{eqnarray*}
  \item [$B_2.$] We now assume that $|\gamma|\leq\min\{\lambda,\ \mu\}.$
Using \eqref{aurelie1.5}
\begin{eqnarray*}
|e^{i\rho v_{_h}(\theta)}-e^{i\rho v_{_h}(\varphi)}|
\leq 2 |I_{h^{\rho}}(\gamma)| + 2|J_{h^{\rho}}(\gamma)|.
\end{eqnarray*}
As in \eqref{bordeaux66},
\begin{equation}
I_{h^{\rho}}(\gamma)\asymp \frac{A_{h^{\rho}}(\gamma)}{h^{2\rho}(\gamma)}.
\end{equation}
Thus
\begin{eqnarray}\label{halwa1}
h^{\rho}(\gamma)|e^{i\rho v_{_h}(\theta)}-e^{i\rho v_{_h}(\varphi)}|
\lesssim  h^{-\rho}(\gamma)A_{h^{\rho}}(\gamma) + h^{\rho}(\gamma)|J_{h^{\rho}}(\gamma)|.
\end{eqnarray}
Using \eqref{condition1} we can compute, as it was done for \eqref{bateau2}, that
\begin{equation}\label{halwa2}
h^{\rho}(\gamma)|J_{h^{\rho}}(\gamma)|\lesssim  \omega(|\gamma|).
\end{equation}
The estimates \eqref{halwa1}  and  \eqref{halwa2}, and the the hypothesis  \eqref{condition2} yield finally to the following estimate
$$h^{\rho}(\gamma)|e^{i\rho v_{_h}(\theta)}-e^{i\rho v_{_h}(\varphi)}|\lesssim  \omega(|\gamma|).$$
\end{itemize}

By combining the above cases $\mathbf{A},$ $B_1$ and $B_2,$ and considering the inequality \eqref{sift} we deduce that $O^{\rho}_h\in\cL(\T),$
and hence $O^{\rho}_h\in\cL.$
\end{proof}

The proof of Theorem \ref{asmahan1} is deduced by joining together the results established in Lemma \ref{linda1}, Lemma \ref{linda2}, Lemma \ref{linda3} and Lemma \ref{construire2}.

\subsection*{\bf Acknowledgements.}
I wish to thank the Abdus Salam International Center for Theoretical Physics
of Trieste in Italy, where a part of this paper was achieved during my stay as a post-doctoral researcher in Mathematics section.

\end{document}